\pdfoutput=1

\documentclass[submission]{FPSAC2021}

\newtheorem{theorem}{Theorem}[section]
\newtheorem{proposition}[theorem]{Proposition}
\newtheorem{lemma}[theorem]{Lemma}

\newtheorem{conjecture}[theorem]{Conjecture}

\usepackage{bm}
\usepackage{dblfloatfix}
\usepackage{mathtools}
\usepackage{enumitem}

\newcommand{\commentout}[1]{}

\newcommand{\euler}[2]{\genfrac{\langle}{\rangle}{0pt}{}{#1}{#2}}

\newcommand{\stirlingcycle}[2]{\genfrac{[}{]}{0pt}{}{#1}{#2}}
\newcommand{\stirlingsubset}[2]{\genfrac{\{}{\}}{0pt}{}{#1}{#2}}

\newcommand{\bp}{{\bm{p}}}
\newcommand{\bq}{{\bm{q}}}
\newcommand{\bw}{{\bm{w}}}
\newcommand{\bfx}{{\mathbf{x}}}

\newcommand{\be}{\begin{equation}}
\newcommand{\ee}{\end{equation}}

\newcommand{\bA}{\bm{A}}
\newcommand{\bS}{\bm{S}}
\newcommand{\bT}{{\bm{T}}}

\newcommand{\wt}[1]{\mathrm{wt}(#1)}
\newcommand{\bP}{{\sf P}}

\newcommand{\calA}{\mathcal{A}}

\newcommand{\calP}{\mathcal{P}}
\newcommand{\calE}{\mathcal{E}}
\newcommand{\Sym}{{\mathfrak{S}}}

\newcommand{\orp}{<_{\pi}}

\newcommand{\ceq}{\coloneqq}

\newcommand{\myge}{\succeq}

\newcommand{\R}{{\mathbb R}}
\newcommand{\Z}{{\mathbb Z}}
\newcommand{\N}{{\mathbb N}}

\newcommand{\tri}{\bigtriangleup}
\newcommand{\triceil}[1]{\lceil #1 \rceil^{\rm tri}}
\newcommand{\tridefect}[1]{\{ #1 \}^{\rm tri}}

\newcommand{\smallest}[1]{\mathop{\rm smallest}\nolimits(#1)}
\newcommand{\largest}[1]{\mathop{\rm largest}\nolimits(#1)}

\newcommand*{\Scale}[2][4]{\scalebox{#1}{$#2$}}

\def\reff#1{(\protect\ref{#1})}

\title{Coefficientwise total positivity of some matrices defined by linear recurrences}

\author[X.~Chen, B.~Deb, A.~Dyachenko, T.~Gilmore, and A.~D.~Sokal]{
		Xi Chen\addressmark{1,2}, 
		 Bishal Deb\addressmark{2}, 
		 Alexander Dyachenko\addressmark{2,3}, 
		 Tomack Gilmore\addressmark{2},
    \and Alan D.~Sokal\addressmark{2,4}}

\address{\addressmark{1}School of Mathematical Sciences, Dalian University of
                   Technology, Dalian 116024, CHINA \\ 
         \addressmark{2}Department of Mathematics, University College London,
                    London WC1E 6BT, UK \\
         \addressmark{3}Department of Mathematics and Statistics,
                    Universit\"at Konstanz, \hbox{D-78457 Konstanz, GERMANY} \\
         \addressmark{4}Department of Physics, New York University,
                    New York, NY 10003, USA}

\received{\today}

\abstract{
   We exhibit a lower-triangular matrix of polynomials $\bm{T}(a,c,d,e,f,g)$
   in~six indeterminates that appears empirically to be
   coefficientwise totally positive, and which includes
   as a special case the Eulerian triangle.
   We prove the coefficientwise total positivity of $\bm{T}(a,c,0,e,0,0)$,
   which includes the reversed Stirling subset triangle.
}

\keywords{Total positivity, coefficientwise total positivity,
   Eulerian triangle.}

\begin{document}

\maketitle

\section{Introduction}

A finite or infinite matrix with integer or real coefficients is
called \emph{totally positive} if all its minors are nonnegative,
and \emph{strictly totally positive} if all its minors are
strictly positive.\footnote{
   {\bf Warning:}  Many authors (e.g.\ \cite{Gantmacher_02,Fomin_00,Fallat_11})
   use the terms "totally nonnegative" and "totally positive"
   for what we have termed "totally positive" and
   "strictly totally positive", respectively.
}
Such matrices have a wide variety of applications across pure and
applied mathematics; background material on this topic can be
found in~\cite{Karlin_68,Gantmacher_02,Pinkus_10,Fallat_11}.
Many interesting lower-triangular matrices (hereafter simply referred
to as \emph{triangles}) that arise in combinatorics have been shown
to be totally positive:
well-known examples include the binomial coefficients $\binom{n}{k}$,
the Stirling cycle numbers $\stirlingcycle{n}{k}$,
and the Stirling subset numbers $\stirlingsubset{n}{k}$.
But there are also many other combinatorially interesting triangles
that appear to be totally positive but for which we have no proof.
Foremost among these is what we call the "clean Eulerian triangle"
\be
        \bA
        \;=\;
        \biggl( \euler{n}{k}^{\! \rm clean} \biggr)_{\!\! n,k \ge 0}
        \;=\;
        \Scale[0.85]{
        \begin{bmatrix}
        1 &    &    &    &    &    &           \\
        1 &  1 &    &    &    &    &           \\
        1 &  4 &  1 &    &    &    &           \\
        1 &  11 &  11 &  1 &    &    &           \\
        1 &  26 &  66 &  26 &  1 &    &           \\
        1 &  57 &  302 &  302 &  57 &  1 &           \\
        \vdots &  \vdots &  \vdots &  \vdots &  \vdots &  \vdots  & \ddots  \\
        \end{bmatrix}
        }
    \,,
\ee
which was conjectured by Brenti~\cite{Brenti_96} to be totally positive,
already a quarter of a century ago.\footnote{
   Note that there exist several different conventions for the
   Eulerian triangle. 
   For our purposes, the "clean" version defined here
   is the most convenient, as it has 1's both on the diagonal
   and in the zeroth column and is reversal-symmetric
   (i.e.\ $\euler{n}{k}^{\! \rm clean} = \euler{n}{n-k}^{\! \rm clean}$).
   It is easy to see that the other versions are totally positive
   if and only if the "clean" one is.
}
Here $\euler{n}{k}^{\! \rm clean}$
is the number of permutations of $[n+1]$ with $k$ excedances (or $k$ descents),
or the number of increasing binary trees on the vertex set $[n+1]$
with $k$ left children. These numbers satisfy the recurrence
\be
	\euler{n}{k}^{\! \rm clean} \;=\; (n-k+1) \,
	\euler{n-1}{k-1}^{\! \rm clean} \:+\: (k+1) \,
	\euler{n-1}{k}^{\! \rm clean} \label{eq.eulerian.recurrence}
\ee
for $n\geq 1$,
with initial condition $\euler{0}{k}^{\! \rm clean} = \delta_{k0}$.

\begin{conjecture}[{\protect\cite[Conjecture~6.10]{Brenti_96}}]
  \label{conj.eulerian}
	The clean Eulerian triangle $\bA$ is totally positive.
\end{conjecture}

A similar problem concerns the reversed Stirling subset triangle. Recall that the Stirling subset number $\stirlingsubset{n}{k}$
is the number of partitions of an $n$-element set into $k$ non-empty blocks
\cite[A048993/A008277]{OEIS}.
We then write $\stirlingsubset{n}{k}^{\rm rev} = \stirlingsubset{n}{n-k}$.
The reversed Stirling subset triangle is \cite[A008278]{OEIS}
\be
   \bS^{\rm rev}
   \;=\;
   \biggl( \stirlingsubset{n}{k}^{\! \rm rev} \biggr)_{\!\! n,k \ge 0}
   \;=\;
   \Scale[0.85]{
   \begin{bmatrix}
	1 &    &    &    &    &    &         \\
	1 &  0 &    &    &    &    &         \\
	1 &  1 &  0 &    &    &    &         \\
	1 &  3 &  1 &  0 &    &    &         \\
	1 &  6 &  7 &  1 &  0 &    &         \\
	1 &  10 &  25 &  15 &  1 &  0 &         \\
	\vdots &  \vdots &  \vdots &  \vdots &  \vdots &  \vdots &   \ddots  \\
   \end{bmatrix}
   }
   \,.
\ee
These numbers satisfy the recurrence
\be
  \stirlingsubset{n}{k}^{\!\rm rev}
   \;=\;
   (n-k) \, \stirlingsubset{n-1}{k-1}^{\!\rm rev}
   \:+\:
   \stirlingsubset{n-1}{k}^{\!\rm rev}
 \label{eq.revstirling.recurrence}
\ee
for $n\geq 1$,
with initial condition $\stirlingsubset{0}{k}^{\!\rm rev} = \delta_{k0}$.
Please note that the total positivity of a lower-triangular matrix
does \emph{not} in general imply the total positivity of its reversal.
Nevertheless we conjecture:

\begin{conjecture}
	\label{conj.revstirling}
The reversed Stirling subset triangle $\bS^{\rm rev}$ is totally positive.
\end{conjecture}

In this extended abstract we present a more general triangle comprised
of \emph{polynomial} entries in \emph{six} indeterminates that
appears empirically to be \emph{coefficientwise} totally positive
and that yields, under suitable specialisations,
both $\bA$ and $\bS^{\rm rev}$.
We do not yet have any proof that this more general triangle
is totally positive; indeed, we do not yet have any proof
of Conjecture~\ref{conj.eulerian}.
But we are able to prove a special case that includes
a generalisation of Conjecture~\ref{conj.revstirling}.

Before stating our main conjecture,
we extend the notion of total positivity to matrices whose elements
are polynomials in one or more indeterminates $\bfx$.
We equip the polynomial ring $\R[\bfx]$ with the
\emph{coefficientwise partial order}:
that is, we say that $P$ is nonnegative (and write $P \myge 0$)
in case $P$ is a polynomial with nonnegative coefficients.
We then say that a matrix with entries in $\R[\bfx]$ is
\emph{coefficientwise totally positive}
if all of its minors are polynomials with nonnegative coefficients.

Comparing recurrences~\eqref{eq.eulerian.recurrence}
and~\eqref{eq.revstirling.recurrence}
invites us to consider the more general linear recurrence
\be
   T(n,k)  \;=\;  [a(n-k) + c] \, T(n-1,k-1) \:+\: (dk + e) \, T(n-1,k) 
 \label{eq.TT2e.b=0}
\ee
for $n \geq 1$, with initial condition $T(0,k) = \delta_{k0}$.
Here $a,c,d,e$ could be integers or real numbers,
but we prefer to treat them as algebraic indeterminates.
Thus, the elements of the matrix $\bT=(T(n,k))_{n,k\geq 0}$
belong to the polynomial ring $\Z[a,c,d,e]$, and we conjecture:

\begin{conjecture}
   \label{conj.TT2e}
The lower-triangular matrix $\bT = \bigl( T(n,k) \bigr)_{n,k \ge 0}$
defined by~\eqref{eq.TT2e.b=0}
is coefficientwise totally positive in the indeterminates $a,c,d,e$.
\end{conjecture}

\noindent
In particular, Conjecture~\ref{conj.eulerian} would follow
by specialising $(a,c,d,e) = (1,1,1,1)$,
while Conjecture~\ref{conj.revstirling} would follow
by specialising $(a,c,d,e) = (1,0,0,1)$.

This, however, is not the end of the story.
Inspired partly by the work of Brenti~\cite{Brenti_95}
and partly by our own experiments,
we were led to consider the \emph{more} general recurrence\hspace*{-2mm}
\be
   T(n,k)  \;=\;  [a(n-k) + c] \, T(n-1,k-1) \:+\:
                  (dk + e) \, T(n-1,k) \:+\:
                  [f(n-2) + g] \, T(n-2,k-1)
 \label{eq.TT5a.intro}
\ee
for $n\geq 1$,
with initial conditions $T(0,k) = \delta_{k0}$ and $T(-1,k) = 0$.
Again, we treat $a,c,d,e,f,g$ as algebraic indeterminates,
so that the matrix elements $T(n,k)$ belong to
the polynomial ring $\Z[a,c,d,e,f,g]$.
Note that this family is invariant under the reversal $k\to n-k$
by interchanging $(a,c)\leftrightarrow(d,e)$ and leaving $f$ and $g$ unchanged:
\be
   T(n,k;\,a,c,d,e,f,g)  \;=\;  T(n,n-k;\,d,e,a,c,f,g)
   \;.
 \label{eq.TT5a.reversal}
\ee
Our main conjecture is the following:

\begin{conjecture}
   \label{conj.TT5a}
The lower-triangular matrix $\bT = \bigl( T(n,k) \bigr)_{n,k \ge 0}$
defined by~\eqref{eq.TT5a.intro}
is coefficientwise totally positive in the indeterminates $a,c,d,e,f,g$.
\end{conjecture}

Unfortunately, for the time being,
Conjectures~\ref{conj.eulerian}, \ref{conj.TT2e} and \ref{conj.TT5a}
remain unproven.
(We have verified Conjecture~\ref{conj.TT5a} up to $13 \times 13$;
this computation took 109 days CPU time.)
The rest of this extended abstract is devoted to proving
the following special case of Conjecture~\ref{conj.TT2e},
which is of some interest in its own right:

\begin{theorem}\label{thm.Tace}
The matrix $\bT=(T(n,k))_{n,k\geq 0}$ specialised to $d=f=g=0$
is coefficientwise totally positive.
\end{theorem}

The triangle that appears in Theorem~\ref{thm.Tace}
is a generalisation of the reversed Stirling subset triangle,
and reduces to it when $(a,c,e) = (1,0,1)$;
this proves Conjecture~\ref{conj.revstirling}.
In what follows we write $\bT(a,c,d,e,f,g)$
for the matrix defined by \eqref{eq.TT5a.intro},
and $\bT(a,c,d,e) = \bT(a,c,d,e,0,0)$
for the matrix defined by \eqref{eq.TT2e.b=0}.

It is possible to prove Theorem~\ref{thm.Tace} in at least two different ways:
one algebraic, the other combinatorial.
In this extended abstract we take the combinatorial path,
leaving the algebraic arguments to a longer paper
(currently under construction).
Section~\ref{sec.comb.interp} establishes combinatorial interpretations
of the entries of $\bT(a,c,0,e)$ and $\bT(0,c,d,e)$
as generating polynomials for set partitions with suitable weights.
In Section~\ref{sec.planarnetworks} we present a planar network $D'$
and show --- by two different arguments ---
that the corresponding path matrix is equal to $\bT(a,c,0,e)$;
Theorem~\ref{thm.Tace} then follows by the Lindstr\"om--Gessel--Viennot lemma.

\section{Set partitions and the matrices $T(a,c,0,e)$ and $T(0,c,d,e)$}
   \label{sec.comb.interp}

From the fundamental recurrence
$\stirlingsubset{n}{k} =
 \stirlingsubset{n-1}{k-1} \,+\, k \, \stirlingsubset{n-1}{k}$
for the Stirling subset numbers
and its consequence \reff{eq.revstirling.recurrence}
for the reversed Stirling subset numbers,
we see that the Stirling and reversed Stirling numbers
correspond to the matrix $\bT(a,c,d,e)$ 
with $(a,c,d,e) = (0,1,1,0)$ and $(1,0,0,1)$, respectively.
Moreover, if one considers instead
$\stirlingsubset{n+1}{k+1}$ and $\stirlingsubset{n+1}{k}^{\rm rev}$,
then these matrices correspond to $\bT(a,c,d,e)$
with $(a,c,d,e) = (0,1,1,1)$ and $(1,1,0,1)$, respectively.
We will now show how to generalise the combinatorial interpretations
of $\stirlingsubset{n+1}{k+1}$ and $\stirlingsubset{n+1}{k}^{\rm rev}$
in terms of set partitions to $\bT(0,c,d,e)$ and $\bT(a,c,0,e)$.

We write $\Pi_n$ (resp.\ $\Pi_{n,k}$)
for the set of all partitions of the set $[n]$ into nonempty blocks
(resp.\ into exactly $k$ nonempty blocks).
For $i \in [n]$ and $\pi \in \Pi_n$,
we write $\smallest{\pi,i}$ for the smallest element
of the block of $\pi$ that contains $i$.
We then have:

\begin{proposition}[Interpretation of $\bT(0,c,d,e)$ and $\bT(a,c,0,e)$
                    in terms of set partitions]
   \label{prop.deb}
\hfill\break
\vspace*{-4mm}
\begin{enumerate}[label=(\roman*)]
   \item  The matrix $\bT = \bT(0,c,d,e)$ has the combinatorial
interpretation
\vspace*{-2mm}
\be
   T(n,k)  \;=\;  \sum_{\pi \in \Pi_{n+1,k+1}} \prod_{i=2}^{n+1} w_\pi(i)
 \label{eq.prop.deb.a1}
\vspace*{-1mm}
\ee
where
\be
   w_\pi(i)
   \;=\;
   \begin{cases}
      e  & \textrm{if $\smallest{\pi,i} = 1$}  \\[-0.5mm]
      c  & \textrm{if $\smallest{\pi,i} = i$}  \\[-0.5mm]
      d  & \textrm{if $\smallest{\pi,i} \neq 1,i$}  \\
   \end{cases}
 \label{eq.prop.deb.a2}
\ee
   \item  The matrix $\bT = \bT(a,c,0,e)$ has the combinatorial
interpretation
\vspace*{-2mm}
\be
   T(n,k)  \;=\;  \sum_{\pi \in \Pi_{n+1,n+1-k}} \prod_{i=2}^{n+1} w_\pi(i)
 \label{eq.prop.deb.b1}
\vspace*{-1mm}
\ee
where
\be
   w_\pi(i)
   \;=\;
   \begin{cases}
      c  & \textrm{if $\smallest{\pi,i} = 1$}  \\[-0.5mm]
      e  & \textrm{if $\smallest{\pi,i} = i$}  \\[-0.5mm]
      a  & \textrm{if $\smallest{\pi,i} \neq 1,i$}  \\
   \end{cases}
 \label{eq.prop.deb.b2}
\ee
\end{enumerate}
\end{proposition}

Please
note that if one restricts
a partition $\pi \in \Pi_{n+1}$ to $[m]$ for some $m < n+1$
--- let us call the result $\pi_m \in \Pi_m$ ---
then $w_\pi(i) = w_{\pi_m}(i)$ for $2 \le i \le m$,
because $\smallest{\pi,i} = \smallest{\pi_m,i}$.
This fact will play a key role in justifying the recurrences.

\begin{proof}[Proof of Proposition~\ref{prop.deb}]
To prove (i)  we will show that the quantities $T(n,k)$ defined by
\reff{eq.prop.deb.a1}/\reff{eq.prop.deb.a2} satisfy the desired recurrence.
Part (ii) follows immediately from (i)
by way of the reversal identity~\eqref{eq.TT5a.reversal} with $f=g=0$.

In a partition $\pi \in \Pi_{n+1,k+1}$,
consider the status of the element $n+1$
and what remains when it is deleted. If $n+1$ is a singleton, then it gets a weight $c$,
and what remains is a partition of $[n]$ with $k$ blocks,
in which each element gets the same weight as it did in $\pi$.
This gives a term $c \, T(n-1,k-1)$. If instead $n+1$ belongs to the block containing 1, then it gets a weight $e$,
and what remains is a partition of $[n]$ with $k+1$ blocks,
in which each element gets the same weight as it did in $\pi$.
This gives a term $e \, T(n-1,k)$. Finally, if $n+1$ belongs to a block whose smallest element lies in
$\{2,3,\ldots,n\}$, then it gets a weight $d$,
and what remains is a partition of $[n]$ with $k+1$ blocks,
in which each element gets the same weight as it did in $\pi$.
There are $k$ blocks not containing~1 to which the element $n+1$
could have been attached.
This gives a term $dk \, T(n-1,k)$.
Summing these terms gives the desired recurrence.
\end{proof}

Here is another recurrence satisfied by these matrices,
which will be useful later:

\begin{lemma}[Alternate recurrences for $\bT(0,c,d,e)$ and $\bT(a,c,0,e)$]
   \label{lemma.deb.recurrence}
\hfill\break
\vspace*{-6mm}
\begin{enumerate}[label=(\roman*)]
   \item  The matrix $\bT = \bT(0,c,d,e)$ satisfies the recurrence
\vspace*{-1mm}
\be
   T(n,k)  \;=\;  e \, T(n-1,k) \:+\: 
                  \sum_{m=0}^{n-1} \binom{n-1}{m} d^m c \, T(n-1-m,k-1)
\vspace*{-1mm}
 \label{eq.lemma.deb.recurrence.1}
\ee
for $n \ge 1$, where $T(n,k) \ceq 0$ if $n<0$ or $k<0$.
   \item  The matrix $\bT = \bT(a,c,0,e)$ satisfies the recurrence
\vspace*{-1mm}
\be
   T(n,k)  \;=\;  c \, T(n-1,k-1) \:+\: 
                  \sum_{m=0}^{n-1} \binom{n-1}{m} a^m e \, T(n-1-m,k-m)
\vspace*{-1mm}
 \label{eq.lemma.deb.recurrence.2}
\ee
for $n \ge 1$, where $T(n,k) \ceq 0$ if $n<0$ or $k<0$.
\end{enumerate}
\end{lemma}

\begin{proof}
(i)  Use the interpretation of Proposition~\ref{prop.deb}(i),
and consider the status of element $n+1$.
If it belongs to the block containing~1, then it gets a weight~$e$,
and what remains is a partition of $[n]$ with $k+1$ blocks;
this gives a term $e \, T(n-1,k)$.
Otherwise, it belongs to a block of size $m+1$ where $0 \le m \le n-1$.
We choose the other $m$ elements of this block in $\binom{n-1}{m}$ ways;
then the smallest element of this block gets weight~$c$,
and the other $m$ elements get weight~$d$.
What remains is a partition of an $(n-m)$-element set with $k$~blocks,
corresponding to $T(n-1-m,k-1)$.

(ii) follows immediately from (i) by the reversal identity.
\end{proof}

We remark that these recurrences, supplemented by the initial
condition $T(0,k) = \delta_{k0}$, completely determine the matrices.


\section{Planar networks and total positivity}\label{sec.planarnetworks}

One very useful tool in proving the total positivity of a matrix
is the famous Lindstr\"om--Gessel--Viennot (LGV) lemma
\cite[Chapter~32]{Aigner_18}.
Consider an acyclic digraph $\cal D$ equipped with edge weights $w_e$
and a distinguished set of sources $U:=\{u_0,u_1,\ldots\}$
and sinks $V:=\{v_0,v_1,\ldots\}$.
The weight $w(\mathcal{P})$ of a path $\mathcal{P}$
is the product of its edge weights;
and we define the \emph{path matrix} $\bm{P} :=(P(u_n\to v_k))_{n,k\geq 0}$ by
$
  P(u_n\to v_k)
  \ceq
  \sum_{\mathcal{P}\colon u_n\to v_k}  w(\mathcal{P})
$.
Now assume further that the digraph $\cal D$ is planar
and that the sources and sinks lie on the boundary of $\cal D$
in the order ``first $U$ in reverse order, then $V$ in order'';
we refer to this setup as a \emph{planar network}.
Then the collection of sources and sinks is \emph{fully compatible}
in the sense that, for any subset of
sources $u_{n_1},\dots,u_{n_r}$ (with $n_1< n_2< \cdots< n_r$)
and sinks $v_{k_1},\dots,v_{k_r}$ (with $k_1<k_2<\cdots<k_r$),
the only permutation $\sigma \in \Sym_r$ mapping
each source $u_{n_i}$ to the sink $v_{k_{\sigma(i)}}$
that gives rise to a nonempty family of nonintersecting paths in $\cal D$
is the identity permutation.
The LGV lemma then implies that every minor of the path matrix $\bm{P}$
is given by a sum over families of nonintersecting paths
between specified subsets of $U$ and $V$,
where each family has weight $\prod w(\mathcal{P}_i)$.
If furthermore every edge weight $w_e$ is a positive real number,
then $\bm{P}$ is totally positive;
and if every edge weight is a polynomial in some indeterminates $\bfx$
with nonnegative real coefficients,
then $\bm{P}$ is coefficientwise totally positive.
This argument goes back to Brenti~\cite{Brenti_95}.

\begin{figure}[t]
\center
\includegraphics[scale=0.4]{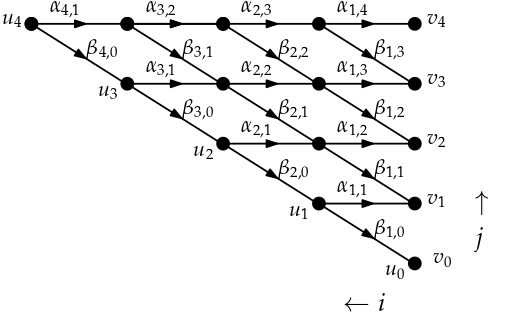} \\[5mm]
\includegraphics[scale=0.45]{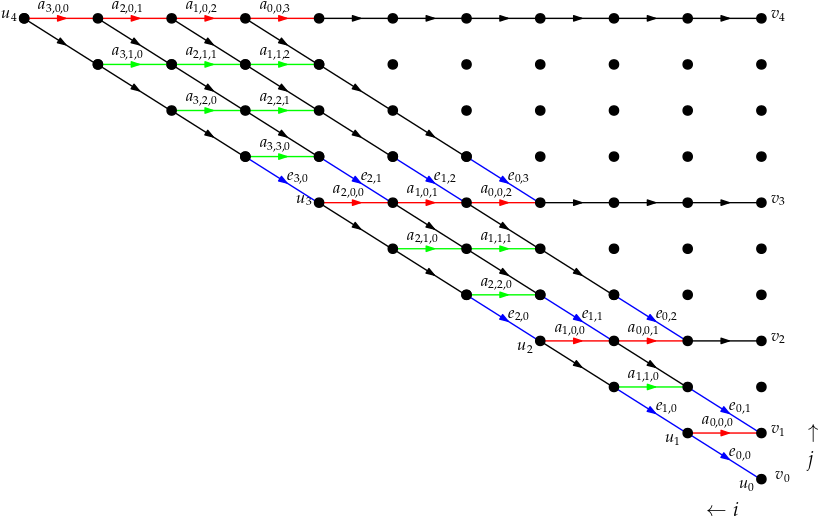}
\caption{
   (a) The standard binomial-like planar network $D$ (above), and
   (b) the planar network $D'$ (below),
   each shown up to source $u_4$ and sink $v_4$.
}
   \label{fig.genericDigraph}
\end{figure}

Figure~\ref{fig.genericDigraph}(a) shows what we call the
\emph{standard binomial-like planar network}, which we denote $D$.
We label the vertices of $D$ by pairs $(i,j)$ with $0 \le i \le j$,
where $i$ increases from right to left
and $j$ increases from bottom to top.
The horizontal directed edge from $(i,j)$ to $(i-1,j)$
[where $1 \le i\leq j$] is given a weight $\alpha_{i,j-i+1}$,
while the diagonal directed edge from $(i,j)$ to $(i-1,j-1)$
[where $1 \le i\leq j$] is given a weight $\beta_{i,j-i}$.
The source vertices are $u_n = (n,n)$
and the sink vertices are $v_k = (0,k)$.

It is easy to see that if the weights are purely $i$-dependent, then 
\be
   P(u_n\to v_k) \;=\;
   \alpha_{n,\bullet} P(u_{n-1}\to v_{k-1}) 
      \,+\, \beta_{n,\bullet} P(u_{n-1}\to v_k)
   \;,
\ee
so that the entries of the corresponding path matrix satisfy a
purely $n$-dependent linear recurrence.
Similarly, if the weights are purely $j$-dependent, then 
\be
   P(u_n\to v_k) \;=\;
   \alpha_{\bullet,k} P(u_{n-1}\to v_{k-1}) 
      \,+\, \beta_{\bullet,k} P(u_{n-1}\to v_{k})
   \;,
\ee
so that the entries of the corresponding path matrix satisfy a
purely $k$-dependent recurrence.
In particular, by setting $\alpha_{i,j}=1$ and $\beta_{i,j}=j$,
we recover a digraph yielding the Stirling subset triangle
$P(u_n\to v_k) = \stirlingsubset{n}{k}$;
and more generally, by setting $\alpha_{i,j}=c$ and $\beta_{i,j}=jd+e$,
we recover $\bT(0,c,d,e)$ and prove its coefficientwise total positivity.
This too goes back to Brenti~\cite{Brenti_95}.

\subsection{The planar network $D'$}

We will now describe a digraph $D'$ that is obtained from $D$
by deleting certain edges (or equivalently, setting their weights to 0),
setting some of the other weights to 1,
and relabelling the remaining weights.
A special role will be played by the \emph{triangular numbers}
$\tri(n) := \binom{n+1}{2}$.
We also define the "triangular ceiling" $\triceil{k}$
to be the smallest triangular number that is $\ge k$,
and the "triangular defect" $\tridefect{k} := \triceil{k} - k$.

For the diagonal edges, we set
\be
   \beta_{i,l}
   \;=\;
   \begin{cases}
       e_{\tri^{-1}(i+l-1)-l,\: l}  &
             \textrm{if $i+l-1$ is triangular and $i+l-1 \ge \tri(l)$} \\
       1  &  \textrm{if $i+l-1$ is not triangular and $i+l-1 \ge \tri(l)$} \\
       0  &  \textrm{in all other cases}
   \end{cases}
\ee
for $i \ge 1$ and $l \ge 0$.
For the horizontal edges, we set
\be
   \alpha_{i,l}
   \;=\;
   \begin{cases}
       a_{\tri^{-1}(\triceil{i+l-1})-l,\: \tridefect{i+l-1},\: l-1}  &
            \textrm{if $\tri^{-1}(\triceil{i+l-1})-l \ge \tridefect{i+l-1}$} \\
       1  &  \textrm{if $i+l-1$ is triangular and $i+l-1 < \tri(l)$} \\
       0  &  \textrm{in all other cases}
   \end{cases}
\ee
for $i,l \ge 1$.
We then delete the edges with zero weight.
Finally, we take the source vertices to be $u_n := (\tri(n),\tri(n))$
and the sink vertices to be $v_{k} := (0,\tri(k))$.
The resulting planar network $D'$ is shown in
Figure~\ref{fig.genericDigraph}(b).


It is clear that every edge of $D'$ either has weight 1
(we call these \emph{black edges})
or else has a unique weight in the set $\calA \cup \calE$, where
$\calA := \{a_{i,j,l} \colon\, (i,j,l)\in\N^3 \textrm{ and } j\leq i\}$
and $\calE := \{e_{i,l} \colon\, (i,l)\in\N^2\}$
(we call these \emph{coloured edges}).
Each path $\calP$ has a weight $w(\calP)$ that is
a monomial in $\Z[\calA,\calE]$.

Let $\bP_{n,k}$ be the set of all paths in $D'$ from $u_n$ to $v_k$.
It is easy to see that $\bP_{n,k}$ is nonempty
if and only if $n \ge k$.
Furthermore, for any two distinct paths $\calP,\calP'$
from $U$ to $V$ in $D'$, we have $w(\calP)\neq w(\calP')$.
Lastly, note that each path $\calP \in \bP_{n,k}$
traverses precisely $n$ coloured edges,
so $w(\calP)$ is a monomial of total degree~$n$.

Applying the Lindstr\"om--Gessel--Viennot lemma to the digraph $D'$,
we can immediately conclude:

\begin{proposition}
   \label{prop.calAE}
The matrix $\bT = (T(n,k))_{n,k \ge 0}$ defined by
$T(n,k) = \sum_{\calP \in \bP_{n,k}} w(\calP)$,
with entries in $\Z[\calA,\calE]$,
is coefficientwise totally positive.
\end{proposition}

The trouble with Proposition~\ref{prop.calAE}
--- as with many applications of Lindstr\"om--Gessel--Viennot ---
is that the set of paths in a digraph can be a rather complicated object;
our goal is to find a simpler combinatorial interpretation.
This can be done either by obtaining a recurrence
that can be compared with Lemma~\ref{lemma.deb.recurrence},
or by constructing an explicit bijection between paths and set partitions.
We shall describe in detail the former approach, and then sketch the latter.

For $0 \le m \le n$, let $u_{n,m} \ceq (\tri(n)-m,\tri(n))$
be the vertex that lies $m$ steps to the right of $u_n$.
We observe that the subnetwork of $D'$ reachable from $u_{n,m}$ is isomorphic
--- after contraction of some black edges,
relabelling $u_n \to u_{n-m}$ and $v_k \to v_{k-m}$ of source and sink vertices,
and relabelling of edge weights ---
to the subnetwork reachable from $u_{n-m}$.
It follows that
\be
   P(u_{n,m} \to v_k)
   \;=\;
   P(u_{n-m} \to v_{k-m}) \bigr|_
    {a_{i,j,l} \to a_{i,j,l+m},\: e_{i,l} \to e_{i,l+m}}
   \;.
 \label{eq.unm.relabeling_identity}
\ee

Now consider a path $\calP$ from $u_n$ to $v_k$.
If the first step is to the right,
we obtain $a_{n-1,0,0}$ times $P(u_{n,1} \to v_k)$.
If the first step is diagonally downwards,
we enter a binomial-like network of size $n-1$,
from which we can emerge on the right wall
at some point $\widehat{u}_{n-1,m} \ceq (\tri(n-1),\tri(n-1)+m)$
for $0 \le m \le n-1$;
from there we follow edges diagonally downwards,
arriving at the point $u_{n-1,m}$
and picking up an extra factor $e_{n-1-m,m}$.
The contribution of the binomial-like network is a bit complicated,
but if we make the specialisation $a_{i,j,l} \to a$ whenever $j > 0$,
then its weight is just $\binom{n-1}{m} a^m$.
We also specialise $a_{i,0,l} \to c_i$ and $e_{i,l} \to e_i$
in order to trivialise the relabellings in \reff{eq.unm.relabeling_identity}.
It follows that with these specialisations the matrix $\bT$
satisfies the recurrence
\be
   T(n,k)
   \;=\;
   c_{n-1} \, T(n-1,k-1)
   \:+\: \sum_{m=0}^{n-1} \binom{n-1}{m} a^m \, e_{n-1-m} \, T(n-1-m,k-m)
 \label{eq.recurrence.acenetwork}
\ee
for $n \ge 1$.
In particular, if $c_i = c$ and $e_i = e$ for all~$i$,
then we recover the recurrence \reff{eq.lemma.deb.recurrence.2}.
Applying Lemma~\ref{lemma.deb.recurrence}(ii), we conclude:

\begin{theorem}
   \label{thm.pathmatrix.ace}
The path matrix $\bT = (T(n,k))_{n,k \ge 0}$ defined by
$T(n,k) = \sum_{\calP \in \bP_{n,k}} w(\calP)$,
with the specialisations $e_{i,l} \to e$,
$a_{i,0,l} \to c$, and $a_{i,j,l} \to a$ for $j > 0$,
coincides with the matrix $\bT(a,c,0,e)$.
\end{theorem}

Combining Proposition~\ref{prop.calAE} with Theorem~\ref{thm.pathmatrix.ace}
proves Theorem~\ref{thm.Tace}.
More generally, Proposition~\ref{prop.calAE} shows that
the matrix $\bT$ defined by the recurrence \reff{eq.recurrence.acenetwork}
is coefficientwise totally positive
in the indeterminates $a$, $(c_i)_{i \ge 0}$ and $(e_i)_{i \ge 0}$.

\subsection{Bijection between paths and set partitions}

We now sketch the bijective approach to proving
Theorem~\ref{thm.pathmatrix.ace},
which is based on a detailed analysis of the paths in the set $\bP_{n,k}$.
The first step is provided by the following lemma,
in which $\wt{\calP}$ denotes the \emph{non-commutative} product of the weights
of $\calP$, taken in the order of traversal.

\begin{lemma}
   \label{lemma.pathword}
Fix integers $n \ge k \ge 0$,
and let $\bw$ be a word in the alphabet $\calA\cup\calE$.
Then~$\bw=\wt{\calP}$ for some path $\calP \in \bP_{n,k}$
if and only if and all of the following conditions hold:
\begin{enumerate}[label=(\roman*)]
\item The first letter of $\bw$ is either $e_{n-1,0}$ or $a_{n-1,j,0}$ where $0\leq j\leq n-1$.
\item The last letter of $\bw$ is either $e_{0,k}$ or $a_{0,0,k-1}$.
\item The letter following $a_{i,j,l}$ is either $e_{i-1,l+1}$ or $a_{i-1,j',l+1}$ where $j\leq j' \leq i-1$.
\item The letter following $e_{i,l}$ is either $e_{i-1,l}$ or $a_{i-1,j,l}$ where $ 0\leq j\leq i-1$.
\end{enumerate}
Furthermore, in this case the word $\bw$ has length $n$
and the path $\cal P$ is unique.
\end{lemma}

\begin{proof}[Sketch of proof]
Parts (i) and (ii) follow from examining the indices of the
first (resp.\ last) coloured edge in $\calP$.
Parts (iii) and (iv) follow from the observation that whenever a
horizontal edge is traversed, $l$ increases by $1$ and $i$ decreases
by $1$; and similarly, when a diagonal coloured edge is traversed,
$l$ remains unchanged and $i$ decreases by $1$. Furthermore, whenever a
coloured horizontal edge is followed immediately by another coloured
edge, the index $j$ is weakly increasing.
\end{proof}

We now construct
a bijection between paths in $D'$
(represented via Lemma~\ref{lemma.pathword} as words)
and set partitions.
%
%
%
Given a set partition $\pi \in \Pi_n$,
we say that an element $i \in [n]$ is\hspace*{-2mm}
\vspace*{-2mm}
\begin{itemize}
   \item an {\em opener}\/ if it is the smallest element of a
      block of size $\ge 2$;  \\[-8mm]
   \item a {\em closer}\/ if it is the largest element of a
      block of size $\ge 2$;  \\[-8mm]
   \item an {\em insider}\/ if it is a non-opener non-closer element of a
      block of size $\ge 3$;  \\[-8mm]
   \item a {\em singleton}\/ if it is the sole element of a block of size 1.
\end{itemize}
\vspace*{-2mm}
Also, for $i \in [n]$ and $\pi \in \Pi_n$,
we write $\mathrm{smallest}(\pi,i)$
for the smallest
element of the block of $\pi$ that contains $i$.

Given a set partition $\pi \in \Pi_{n+1,k}$
consisting of blocks $B_1,\ldots,B_k$,
we define a total order $\orp$ on $[n+1]$
by the following procedure:
start by taking the block containing~1 (we call it~$B_1$) together with the
largest elements of all the other blocks, and put them in increasing order;
then insert all the remaining elements of each block
(other than~$B_1$) in increasing order
immediately preceding its largest element.
For example,
for $\pi = \{ \{1,5,8\}, \{2,3,9\}, \{4,7\}, \{6\} \} \in \Pi_{9,4}$,
the order is $156478239$.
%
%
%
%

Under this total order, $1$ is the smallest element and $n+1$ is the largest;
it can therefore be written as $1 p_1 p_2 \cdots p_n$ where $p_n = n+1$.
%
%
We then define the word associated to a set partition $\pi\in\Pi_{n+1,n-k+1}$
to be $W(\pi) \ceq w_n \cdots w_1$
where
\be
   w_i
   \;\ceq\;
   \begin{cases}
      e_{i-1,\: l_i}       & \textrm{if $\smallest{\pi,p_i} = p_i$} \\
      a_{i-1,\: 0,\: l_i}  & \textrm{if $p_i\in B_1$ \ 
             [i.e.\ $\smallest{\pi,p_i} = 1$]}  \\
      a_{i-1,\: j,\: l_i}  & \textrm{if $\smallest{\pi,p_i}\neq 1,p_i$
                                     and $\largest{\bp,i}_j = p_{i-1}$}
   \end{cases}
 \label{def.wi}
\ee
Here $\largest{\bp,i}_j$ denotes the $j$th largest element
of the set $[2,n+1]\setminus\{p_i,\ldots,p_n\}$.
The index $l_i$ is defined recursively:
we set $l_n=0$, and for $i<n$ we define
$l_{i-1} = l_i$ if $\smallest{\pi,p_i} = p_i$,
and $l_{i-1} = l_i + 1$ otherwise.

\begin{lemma}
   \label{lemma.setpartword}
Given $\pi\in \Pi_{n+1,n-k+1}$, the word $W(\pi)$
consists of letters from $\calA\cup \calE$ and satisfies the conditions
in Lemma~\ref{lemma.pathword}, thereby corresponding to a path
$\calP \in \bP_{n,k}$.
\end{lemma}

\begin{proof}[Outline of proof]
We first verify that $w_i\in \calA\cup \calE$ (this is easy);
then we check the four conditions of Lemma~\ref{lemma.pathword}
(this requires some consideration of cases).
\end{proof}

Lemmas~\ref{lemma.pathword} and \ref{lemma.setpartword} together define
a map $\Phi_{n,k} \colon\, \Pi_{n+1,n-k+1}\to \bP_{n,k}$.

\begin{theorem}
   \label{thm.bijection}
The map $\Phi_{n,k}$ is a bijection of $\Pi_{n+1,n-k+1}$ onto $\bP_{n,k}$.
\end{theorem}

\begin{proof}[Outline of proof]
Given a word $\bw=w_n\cdots w_1$
satisfying the conditions in Lemma~\ref{lemma.pathword},
we construct a set partition $\pi \in \Pi_{n+1,n-k+1}$
satisfying $W(\pi) = \bw$; this will show surjectivity.
We also show that $\pi$ is the unique set partition
with this property, showing injectivity.
We build up $\pi$ by inserting elements into its blocks,
one at a time, as we read the word $\bw$ from left to right,
beginning from $\pi_0 = \{\{1\}\}$ and ending with $\pi_n = \pi$.
Each block will be built up in decreasing order,
starting with its largest element;
indeed, each block other than $B_1$
will be built from start to finish in successive stages of the algorithm.
Whenever we insert an element $q_i \in [2,n+1]$ into a block $B \neq B_1$,
we also declare whether that block is
\emph{finished} (i.e.\ $q_i$ is an opener or a singleton in $\pi$)
or \emph{unfinished} (i.e.\ $q_i$ is a closer or an insider in $\pi$).
At each stage there will be at most one unfinished block.
We show \emph{a~posteriori} that $q_i$ equals
the $p_i$ associated to the total order $\orp$.

When we read a letter $w_i$,
we choose an element $q_i \in [2,n+1]$
that is not already contained in $\pi_{n-i}$,
and insert it into $\pi_{n-i}$ in one of five ways:
insert $q_i$ into the block $B_1$;
insert $q_i$ as an opener into an unfinished block $B\neq B_1$;
insert $q_i$ as an insider into an unfinished block $B\neq B_1$;
create a new block containing $q_i$ as a singleton;
or create a new block containing $q_i$ as a closer.
The result is called $\pi_{n-i+1}$.

To construct the sequence $\bq = q_n \cdots q_1$,
we start from $q_n = n+1$.
Then, for $i < n$, we proceed inductively:
if $w_{i+1} = a_{i,j,l}$ with $j > 0$,
we set $q_i$ to be the $j$th largest element
of the set $[2,n+1] \setminus \{q_{i+1},\ldots,q_n\}$;
otherwise we set $q_i$ to be the largest element
of $[2,n+1] \setminus \{q_{i+1},\ldots,q_n\}$.

The elements $q_n,\ldots,q_1$ are inserted successively
into the set partition as follows:
By Lemma~\ref{lemma.pathword}
there are three possibilities for the letter $w_i$:
$e_{i-1,l}$, $a_{i-1,0,l}$, or $a_{i-1,j,l}$ for some $1\leq j\leq i-1$.

Case 1: $w_i = e_{i-1,l}$.
If there is an unfinished block, we insert $q_i$ into that block as an opener;
otherwise, we create a new block with $q_i$ as a singleton.

Case 2: $w_i = a_{i-1,0,l}$.
We insert $q_i$ into block $B_1$.

Case 3: $w_i = a_{i-1,j,l}$ for some $1\leq j\leq i-1$.
If there is an unfinished block, we insert $q_i$ into that block as an insider;
otherwise, we create a new block with $q_i$ as a closer.

We then prove: 1) the claims about the order in which the blocks are built;
2) that $\bp = \bq$; 3) that $W(\pi) = \bw$; and 4) that the map is injective.
All these steps require some consideration of cases.
\end{proof}

\begin{proof}[Second proof of Theorem~\ref{thm.pathmatrix.ace}]
The definition \reff{def.wi} tells us that,
within each word $\bw$,
the letters $a_{i,0,l}$ correspond to elements in $B_1$,
and the letters $e_{i,l}$ (resp.\ $a_{i,j,l}$ for $j>0$)
correspond to minimal (resp.\ non-minimal) elements of blocks $B \neq B_1$.
After the specialisations $e_{i,l} \to e$,
$a_{i,0,l} \to c$, and $a_{i,j,l} \to a$ for $j > 0$,
by Proposition~\ref{prop.deb}(ii) this is precisely the matrix $\bT(a,c,0,e)$.
\end{proof}

\section*{Acknowledgements}
We wish to thank Sergey Fomin for helpful correspondence.
This research was supported in part by
the U.K.~Engineering and Physical Sciences Research Council grant EP/N025636/1,
a fellowship from the China Scholarship Council,
and a fellowship from the Deutsche Forschungsgemeinschaft.



\end{document}